\documentclass[a4paper,10pt]{article}
\usepackage{url}
\usepackage{latexsym}
\usepackage{amssymb}
\usepackage{amsmath}
\usepackage{amsthm}
\usepackage{times}
\usepackage{xspace}
\usepackage[dvips]{graphicx}
\usepackage{epsf}
\usepackage{amsfonts}
\usepackage{amssymb}

\newtheorem{lem}{Lemma}
\def\qed{\hfill \ifhmode\unskip\nobreak\fi\quad\ifmmode\Box\else$\Box$\fi\\ }
{\bfseries}{\itshape}
\newtheorem{cor}{Corollary}{\bfseries}{\itshape}
{\bfseries}{\itshape}
{\bfseries}{\itshape}
\newtheorem{theorem}{Theorem}{\bfseries}{\itshape}
\newtheorem{lemma}{Lemma}{\bfseries}{\itshape}
\newtheorem{definition}{Definition}{\bfseries}{\itshape}

\title{Hadwiger Number and the Cartesian Product  of Graphs}
\author{L. Sunil Chandran
\thanks{Dept. of Computer Science and Automation, Indian Institute of Science,
Bangalore, 560012. Email : {\tt sunil@csa.iisc.ernet.in} Alternate
Email: {sunil.cl@gmail.com} } \and Alexandr Kostochka \thanks
{University of Illinois, Urbana, IL-61801, USA and Sobolev Institute of
Mathematics, Novosibirsk, 630090, Russia. Email: {\tt
kostochk@math.uiuc.edu}. The work of this author is partially supported
by the NSF grant DMS-06-50784.} \and J. Krishnam Raju \thanks{David R.
Cheriton School of Computer Science, University of Waterloo, Email:{\tt
krjampan@cs.uwaterloo.ca}}}
\date{}

\begin{document}
\pagestyle{plain}
\pagenumbering{arabic}
\maketitle

\newcommand{\mr}{\ensuremath{\mathbf{\eta}}}
\newcommand{\len}{\ensuremath{t}}
\newcommand{\minor}{\ensuremath{\preceq}}
\newcommand{\node}[1]{\ensuremath{\langle #1 \rangle }}
\newcommand{\floor}[1]{\ensuremath{\lfloor #1 \rfloor }}
\newcommand{\ceil}[1]{\ensuremath{\lceil #1 \rceil }}
\newcommand{\cart}{\ensuremath{\,\Box\,}}
\newcommand{\brac}[1]{\ensuremath{\lbrace #1 \rbrace }}
\newcommand{\ignore}[1]{}

\begin{abstract}
The Hadwiger number $\mr(G)$ of a graph $G$ is  the largest
integer $n$ for which the complete graph $K_n$ on $n$ vertices is a
minor of $G$. Hadwiger conjectured that for every graph $G$, $\mr(G)
\ge \chi(G)$, where $\chi(G)$ is the chromatic number of $G$. In
this paper, we study  the Hadwiger number of the
Cartesian product $G\cart H$  of graphs.

As the main result of this paper, we prove that
 $\mr(G_1 \cart G_2) \ge h\sqrt{l}\left (1 - o(1) \right )$
for any
two graphs $G_1$ and $G_2$ with $\mr(G_1) = h$ and
$\mr(G_2) = l$.
We show that the above lower bound is asymptotically best
possible. This asymptotically settles a question of Z. Miller
(1978).

As consequences of our main result, we show the following:
\begin{enumerate}
\item Let $G$ be a connected graph. Let the (unique) prime factorization of $G$
be given by $G_1 \cart G_2 \cart ... \cart G_k$. Then
$G$ satisfies Hadwiger's conjecture if $k \ge 2\log{\log{\chi(G)}} + c'$, where
$c'$ is a constant. This improves the $2\log{\chi(G)}+3$ bound  in
\cite{sunns1}.

\item Let $G_1$ and $G_2$ be two graphs such that $\chi(G_1) \ge \chi(G_2) \ge
 c{{\log}^{1.5}}(\chi(G_1))$, where $c$ is a constant. Then $G_1 \cart G_2$
satisfies Hadwiger's conjecture.

\item  Hadwiger's conjecture is true for $G^d$
(Cartesian product of $G$ taken $d$ times) for
every graph $G$ and every $d \ge 2$. This settles a question
by Chandran and Sivadasan
\cite{sunns1}. (They had shown that the Hadiwger's conjecture
is true for $G^d$  if $d \ge 3$.)
\end{enumerate}
\end{abstract}


\section{Introduction}

\subsection{General definitions and notation}
In this paper we only consider undirected simple graphs i.e., graphs without
multiple edges and loops. For a graph $G$, we use $V(G)$ to
denote its vertex set and $E(G)$ to denote its edge set.

A $k$-{\em coloring} of a graph $G(V,E)$ is a function $f:V\rightarrow
\brac{1,2,...,k}$. A $k$-coloring $f$ is {\em proper} if for all edges $(x,y)$ in $G$,
$f(x) \not= f(y)$. A graph is $k$-colorable if it has a proper $k$-coloring.
The {\em chromatic number} $\chi(G)$ is the least $k$ such that $G$ is $k$-colorable.

Let $S_1 ,S_2 \subset V(G)$, such that $S_1 \not= \emptyset$, $S_2 \not=
\emptyset$ and
$S_1 \cap S_2 = \emptyset$. We say that $S_1$ and $S_2$ are adjacent in $G$ if
and only if there exists an edge $(u,v) \in E(G)$ such that $u \in S_1$ and
$v \in S_2$. The edge $(u,v)$ is said to connect $S_1$ and $S_2$.

{\em Contraction} of an edge $e = (x,y)$
is the replacement of  the vertices $x$ and $y$ with a new vertex $z$,
whose incident edges are the edges other than $e$ that were incident to $x$
or $y$. The resulting graph denoted by $G.e$ may be a multigraph, but since
we are only interested in simple graphs, we discard any parallel edges.

A $minor$ $M$ of $G(V,E)$ is  a graph obtained from $G$  by a sequence of
contractions of edges and deletions of edges and vertices.
We call $M$ a minor of $G$ and write $M \minor G$.

It is not difficult to verify that $M \minor G$ if and only if
for each vertex $x \in V(M)$, there exist a set $V_x \subseteq
 V(G)$, such that (1) every  $V_x$ induces a connected subgraph of $G$,
(2) all $V_x$ are disjoint, and (3) for each $(x,y) \in E(M)$,
 $V_x$ is adjacent to $V_y$ in $G$.

The Hadwiger number $\mr(G)$ is the largest integer $h$ such that the complete
graph on $h$ vertices $K_h$ is a minor of $G$. Since every graph on at most $h$ vertices
is a minor of $K_h$, it is easy to see that $\mr(G)$ is the largest integer such
that any graph on at most $\mr(G)$ vertices is a minor of $G$. Hadwiger conjectured
the following in 1943.
\\

\noindent
\textbf{Conjecture:(Hadwiger \cite{Hadwiger43})}
\textit{For every graph $G$, $\mr(G) \ge \chi(G)$, where $\chi(G)$ is the chromatic number of $G$.}
\\

In other words, Hadwiger's conjecture states that if $\mr(G) \le k$, then $G$
is $k$-colorable. It is known to hold for small $k$. Graphs of
Hadwiger number at most 2 are the forests. By a theorem of Dirac \cite{Dirac1},
the graphs with Hadwiger number at most 3 are the series-parallel graphs. Graphs
with Hadwiger number at most 4 are characterized by Wagner \cite{Wagner2}.
The case $k$ = 4 of the conjecture 
implies the Four Color Theorem because any planar graph has no $K_5$ minor. On the
other hand, Hadwiger's conjecture for the case $k$ = 4 follows from the four
color theorem and a structure theorem of Wagner \cite{Wagner2}.
 Hadwiger's conjecture for
$k=5$ was settled by Robertson et al. \cite{RobertsonSeymourThomas93}. The case
$k=6$ onwards is still open.

Since Hadwiger's conjecture in the general case is still open, researchers
have shown interest to derive lower bounds for Hadwiger number in terms of the
chromatic number.
Mader\cite{Mad2} showed (improving an earlier
result of Wagner \cite{Wagner1}) that for any graph $G$,
$\mr(G) \ge \frac{\chi(G)}{16\log(\chi(G))}$. Later,
 Kostochka \cite{Kostochka82} and  Thomason~\cite{AGThoma}
 independently  showed that there exists a constant $c$, such that for
any graph $G$, $\mr(G) \ge \frac{\chi(G)}{c\cdot\sqrt{\log(\chi(G))}}$.

It is also known that Hadwiger's
conjecture is true for almost all graphs on $n$ vertices.

Improving on previous results by other authors,
K\"uhn and Osthus \cite{KuhOs} showed
that if the girth (i.e., the length of a shortest cycle)
is at least $g$ for some odd $g$ and the minimum degree $\delta$
is at least 3, then $\mr(G) \ge \frac{{c(\delta)}^{(g+1)/4}}{\sqrt{\log{\delta}}}$.
As a
consequence of this result, K\"uhn and Osthus\cite{KuhOs} showed that Hadwiger's
conjecture is
true for $C_4$-free graphs of sufficiently large chromatic number. (Here $C_4$
denotes a cycle of length $4$)

\subsection{ The Cartesian product of graphs}

Let $G_1$ and $G_2$ be two undirected graphs,
where the vertex set of $G_1$ is
$\brac{0,1,\cdots,n_1-1}$ and the vertex set of $G_2$ is $\brac{0,1,\cdots$,$n_2-1}$.
The {\em Cartesian  product}, $G_1 \cart G_2$, of $G_1$ and $G_2$ is a graph
with the vertex set $V = \brac{0,1,...,n_1-1} \times
\brac{0,1,...,n_2-1}$ and the edge set  defined as follows.
There is an edge between vertices
$\node{i,j}$
and $\node{i^\prime,j^\prime}$ of $V$ if and only if, either $j = j^\prime$
and $(i,i^\prime) \in E(G_1)$, or $i=i^\prime$ and
$(j,j^\prime)\in E(G_2)$.

In other words, graph products can be viewed in the following way: let the
vertices of $G_1 \cart G_2$ be partitioned into $n_2$ classes $W_1,...,W_{n_2}$,
where  $W_j = \brac{\node{1,j},\cdots,\node{n_1,j}}$ induces a graph that is
isomorphic to $G_1$, where the vertex $\node{i,j}$ corresponds to vertex $i$ of $G_1$.
If edge $(j,j^\prime)$ belongs to $G_2$ then the edges between classes $W_j$
and $W_{j^\prime}$ form a matching such that the corresponding vertices, i.e.,
$\node{i,j}$ and $\node{i,j^\prime}$, are matched. If edge $(j,j^\prime)$ is not present
in $G_2$ then there is no edge between $W_j$ and $W_{j^\prime}$.

It is easy to verify that the Cartesian product
is a commutative and associative operation on graphs.
Due to the associativity, the product of
graphs $G_1,...,G_k$ can be simply written as $G_1 \cart ... \cart G_k$ and
has the following interpretation. If the vertex set of graph $G_i$ is
$V_i$ = $\brac{1,...,n_i}$, then $G_1 \cart ... \cart G_k$ has the vertex set
$V$ = $V_1 \times V_2 \times ... \times V_k$. There is an edge between vertex
$\node{i_1,...,i_k}$
and vertex $\node{i'_1,...,i'_k}$ of $V$ if and only if there is
a position $t,1\le t \le k$, such that $i_1$ = $i'_1$, $i_2$ =
$i'_2$
,...,$i_{t-1}$ = $i'_{t-1}$, $i_{t+1}$ = $i'_{t+1}$,...,
$i_k$ = $i'_k$, and the edge $(i_t, i'_t)$ belongs to graph $G_t$.

We denote the product of graph $G$ taken $k$ times as $G^k$. It is easy to
verify that if $G$ has $n$ vertices and $m$ edges, then $G^k$ has $n^k$ vertices
and $mk\cdot n^{k-1}$ edges. \\

Well known examples of Cartesian products of graphs are the $d$-dimensional
hypercube $Q_d$, which is
isomorphic to $K^d_2$,  and a $d$-dimensional grid, which is isomorphic to
$P_n^d$, where $P_n$ is a simple path on $n$ vertices.

\noindent
\textbf{Unique Prime Factorization(UPF) of graphs: }
A graph $P$ is \textit{prime} with
respect to the Cartesian product operation if and only if $P$ has at least
two vertices and
it is not isomorphic to the product of two non-identity graphs, where an
identity graph is the graph on a single vertex and having no edge. It is
well-known that every connected undirected graph $G$ with at least two vertices has a UPF with respect to
Cartesian product in the sense that if $G$ is not prime then it can be expressed
in a unique way as a product of prime graphs(\cite{Sandy}).
If $G$ can be expressed as the product
$G_1 \cart G_2 \cart ... \cart G_k$, where each $G_i$ is prime, then we say that
the product dimension of $G$ is $k$. The UPF of a given connected graph $G$ can
be found in $O(m\log(n))$ time, where $m$ and $n$ are the number of edges and
number of vertices of $G$ respectively \cite{auren}.\\

Imrich and Klav\^{z}ar have published a book \cite{Sandy}, dedicated
exclusively to the study of graph products. Readers who are interested
to get an introduction to the wealth of profound and beautiful results
on graph products are referred to this book.

We will use the following result by Sabidussi
\cite{Sabi57} (which was rediscovered several times).

\begin{lemma}
\label{max}
$\chi(G_1 \cart G_2) = \max\brac{\chi(G_1), \chi(G_2)}$.
\end{lemma}

\subsection{Our results}

The question of studying the Hadwiger number with respect to the Cartesian
product operation was suggested by Miller in the open problems section of
a 1978 paper \cite{Miller1}. He mentioned a couple of special cases
(such as $\mr(C_n \cart K_2) = 4$  and  $\mr(T \cart K_n) = n+1$, where $C_n$
and $T$ denote a cycle and a tree respectively) and left the general case open.
In this paper, we answer this question asymptotically.
We give the following results.
\vspace{0.32cm} \\
\noindent
\textbf{Result 1. }
Let $G_1$ and $G_2$ be two graphs with $\mr(G_1) = k_1$ and $\mr(G_2) = k_2$.
Then $\mr(G_1 \cart G_2) \ge k_1 \sqrt{k_2}\left ( 1 -o(1) \right )$.
 (Since
the Cartesian product is commutative, we can assume without loss of generality
that $k_1 \ge k_2$). We demonstrate that this lower bound is asymptotically
best possible.

We also show that in general, $\mr(G_1 \cart G_2)$ does not have any upper bound
that depends only on $\mr(G_1)$ and $\mr(G_2)$, by demonstrating graphs $G_1$ and
$G_2$ such that $\mr(G_1)$ and $\mr(G_2)$ are bounded, whereas $\mr(G_1 \cart
G_2)$ grows with the number of vertices of $G_1 \cart G_2$.
\vspace{0.32cm} \\
\noindent
\textit{Remark. }
Note that if the average degrees of $G_1$ and $G_2$ are $d_1$ and
$d_2$ respectively, then the average degree of $G_1 \cart G_2$ is $d_1 + d_2$.
In comparison, by Result 1, the Hadwiger number of $G_1 \cart G_2$ grows much
faster.

Hadwiger's conjecture for Cartesian products of graphs was studied in
\cite{sunns1}. There it was shown that if the product dimension (number of
factors in the unique prime factorization of $G$) is $k$, then Hadwiger's
conjecture is true for $G$ if $k \ge 2\log{\chi(G)}+3$. As a consequence of
Result 1, we are able to improve this bound. We show the following.
\vspace{0.32cm} \\
\noindent
\textbf{Result 2. }
\noindent
Let the (unique) prime factorization of $G$ be $G = G_1 \cart G_2 \cart
\cdots \cart G_k$. Then Hadwiger's conjecture is true for $G$ if $k \ge
2\log(\log(\chi(G))) + c'$, where $c'$ is a constant.

Another consequence of Result 1 is that if $G_1$ and $G_2$ are two graphs
such that $\chi(G_2)$ is not ``too low" compared to $\chi(G_1)$, then
Hadwiger's conjecture is true for $G_1 \cart G_2$. More precisely:
\vspace{0.32cm} \\
\noindent
\textbf{Result 3. }
If $\chi(G_2) \ge c{{\log}^{1.5}}(\chi(G_1))$, where $c$ is a constant,
then Hadwiger's conjecture is true for $G_1 \cart G_2$.

It is easy to see that Result 3 implies the following: Let $G_1$ and $G_2$ be
two graphs such that $\chi(G_1) = \chi(G_2)$. (For example, as in the case
$G_1 = G_2$). Then $G_1 \cart G_2$ satisfies Hadwiger's conjecture if $\chi(G_1)
= \chi(G_2) = t$ is sufficiently large. ($t$ has to be sufficiently large, because
of the constant $c$ involved in Result 3). For this special case, namely
$\chi(G_1) = \chi(G_2)$, we give a different proof (which does not
depend on Result 1), to show that Hadwiger's conjecture is true for $G_1 \cart G_2$.
This proof does not require that $\chi(G_1)$ be sufficiently large.
\vspace{0.32cm} \\
\noindent
\textbf{Result 4. }
Let $G_1$ and $G_2$ be any two graphs such that $\chi(G_1) = \chi(G_2)$. Then
Hadwiger's conjecture is true for $G_1 \cart G_2$.

It was shown in \cite{sunns1} that Hadwiger's conjecture is true for
$G^d$, where $d \ge 3$, for any graph $G$. As a consequence of Result 4, we are
able to sharpen this result.
\vspace{0.32cm} \\
\noindent
\textbf{Result 5. } For any graph $G$ and every $d \ge 2$,
Hadwiger's Conjecture is true for $G^d$.

Another author who studied the minors of the Cartesian product of graphs is Kotlov
\cite{Kotlov}. He showed that for every bipartite graph $G$, the strong product
(\cite{Sandy})
$G \boxtimes K_2$ is a minor of $G \cart C_4$. ($K_2$ and $C_4$ are an edge and a
$4$-cycle respectively). As a consequence of this he showed that $\mr(K_2^d)
\ge 2^{\frac{d+1}{2}}$.

\section{Hadwiger Number for $G_1 \cart G_2$}

\subsection{Lower bound on $\mr(G_1 \cart G_2)$}

\noindent
The following Lemma is not difficult to prove.
\begin{lemma}
\cite{sunns1}
\label{minor}
If $M_1 \preceq G_1$ and $M_2 \preceq G_2$, then
$M_1 \cart M_2 \preceq G_1 \cart G_2$.
\end{lemma}

\noindent
Let $G_1$ and $G_2$ be two graphs such that $\mr(G_1) = h$ and
$\mr(G_2) = l$, with $h \ge l$. In this section we show that
$\mr(G_1 \cart G_2) \ge h \sqrt{l}\left (1-o(1) \right )$.
Since by Lemma \ref{minor}, $K_h \cart K_l \minor G_1 \cart G_2$ it is
sufficient to prove that $\mr(K_h \cart K_l) \ge h
\sqrt{l}\left (1 - o(1) \right )$.


\begin {definition}
An {\em affine plane ${\cal A}$ of order $m$} is a family
$\{A_{i,t}\;:\; i=1,\ldots,m+1,\; t=1,\ldots,m\}$ of $m$-elements
subsets of an $m^2$-element set $A$ such that\\ $$|A_{i,t}\cap
A_{i',t'}|=\left\{
\begin{array}{rl}
1& \mbox{if $i'\neq i$};\\ 0& \mbox{if $i'=i$ and $t'\neq t$.}
\end{array}\right.$$

The sets $A_{it}$ are {\em the lines of ${\cal A}$}. By
definition, for each $i$, the sets $A_{i,1},A_{i,2},\ldots, A_{i,m}$
are disjoint and form a
 partition of $A$.

\end {definition}

 The following fact is widely known (see, e.g.,~\cite{R}):

 \begin{lem}\label{l1} If $m$ is a prime power, then there exists
 an affine plane of order $m$.
 \end{lem}

 It is also known that the set of prime numbers is quite dense.
 The following is a weakening of the result of Iwaniec and Pintz~\cite{IP}.

\begin{lem}\label{l2} For every sufficiently large positive $x$,
the interval $[x,x+x^{0.6}]$  contains a prime number.
 \end{lem}

\begin{cor}\label{cor1} For every sufficiently large positive $x$,
the interval $[x-6x^{0.9},x]$
 contains a number of the form $(p(p+1))^2$, where $p$ is some prime.
 \end{cor}

\begin {theorem}

\label {result2}

 $\eta(K_h \Box K_l) \ge  h \sqrt l (1 - o(1))$.

\end {theorem}

\begin {proof}

Consider $G=K_h\Box K_l$, where $h\geq l$ and $l$ is large.
Let $p$ be the maximum prime
such that $l\geq (p(p+1))^2$.
We view $K_h\Box K_l$ as a set of $h$ copies
of $K_l$. Suppose that $s(p-1)(2p+1)/2\leq h<(s+1)(p-1)(2p+1)/2$.
We neglect some $h-s(p-1)(2p+1)/2$ copies of $K_l$ and partition the
remaining $s(p-1)(2p+1)/2$ copies into $s$ large groups of the same size,
and each of these groups into $(p-1)/2$ groups of size $2p+1$.
In other words, we consider
${\cal S}=\{ K_l(i,j,m)\,:\,i=1,\ldots,s,\; j=1,\ldots,(p-1)/2,\;
m=1,\ldots,2p+1\}$, where
each $K_l(i,j,m)$ is a copy of $K_l$.
For $i=1,\ldots,s,$ let
${\cal S}_i=\{ K_l(i,j,m)\,:\, j=1,\ldots,(p-1)/2;\;
m=1,\ldots,2p+1\}$.

 In  ${\cal S}_1$, we will find
$p^2(p-1)(2p+1)/2$ disjoint sets $M(1,j,m,t)$ ($j=1,\ldots,(p-1)/2$,
$m=1,\ldots,2p+1$, $t=1,\ldots,p^2$)
of size $(p+1)^2$. These sets will have the property that\\
(a) the subgraph $G(M(1,j,m,t))$ induced by $M(1,j,m,t)$ is
connected; \\
(b)  for any two
quadruples $(1,j,m,t)$ and $(1,j',m',t')$, there is a vertex $v$ in $K_l$
such that each of $M(1,j,m,t)$ and $M(1,j',m',t')$ contains a copy of $v$
(in different copies of $K_l$, since our sets are disjoint).

If we manage this, then copying these sets for every $i=2,\ldots,s$,
by (b), we will create $s\ p^2(p-1)(2p+1)/2=p^2h(1-o(1))$ disjoint sets
$M(i,j,m,t)$ that satisfy\\
(a') the subgraph $G(M(i,j,m,t))$ induced by $M(i,j,m,t)$ is
connected; \\
(b')  for any two
quadruples $(i,j,m,t)$ and $(i',j',m',t')$, there is a vertex $v$ in $K_l$
such that each of $M(i,j,m,t)$ and $M(i',j',m',t')$ contains a copy of $v$

So, we go after (a) and (b).

To achieve this, we view the set of vertices of each $K_l$ as the
disjoint union of a "big" square $Q_0$ of size $p^2\times p^2$ with $2p+1$
"small" squares $Q_k$, $k=1,\ldots,2p+1$ of size $p\times p$ and
the reminder $R$ (of size $l-p^4-(2p+1)p^2=l-(p(p+1))^2$) (see
Fig. 1).

\begin{figure}[htbp]
{
  \includegraphics[scale=0.7]{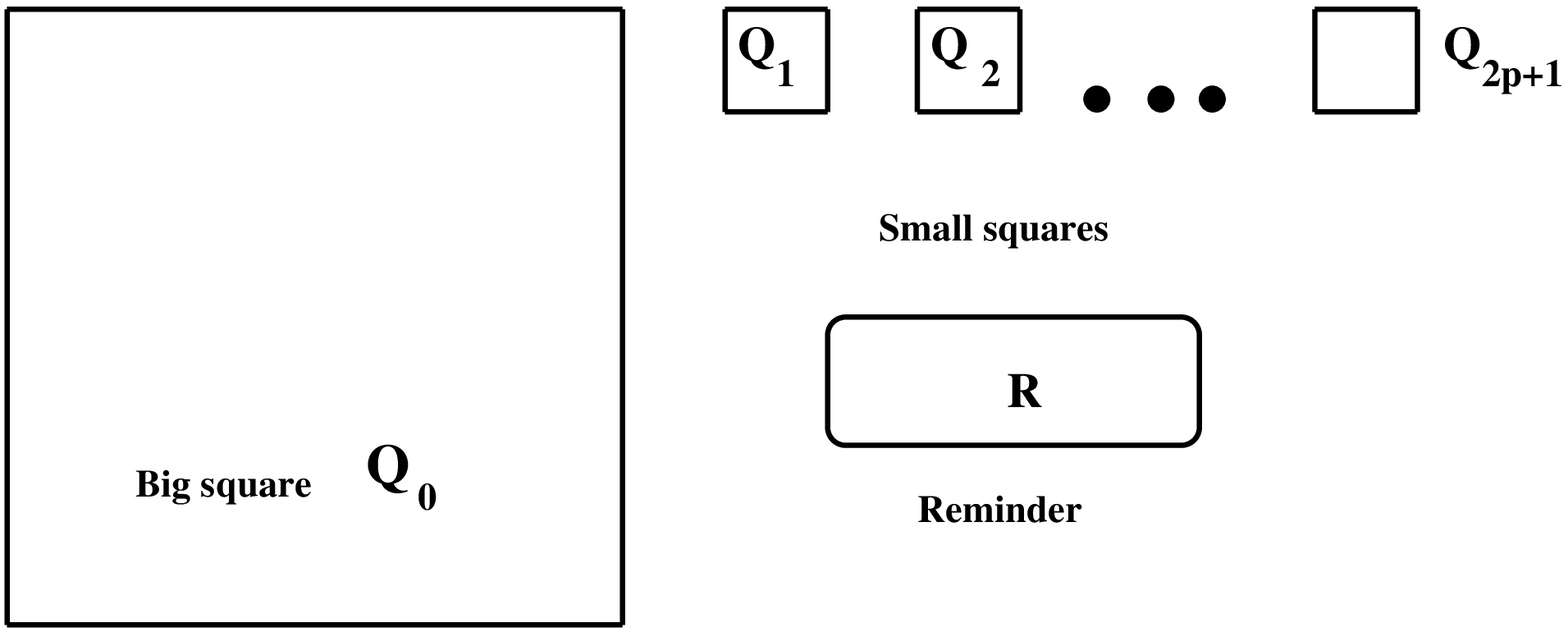}
 }
\end{figure}

By Lemma~\ref{l1}, there exists an affine plane of order $p^2$. We
consider the lines $A_{i,t}$ of this plane as subsets of the big
square $Q_0$. For $i=1,\ldots,(p-1)(2p+1)/2$, we view
$\{A_{i,1},A_{i,2},\ldots,A_{i,p^2}\}$ as a partition of the copy
$Q_0(i)$ of the big square $Q_0$. If $i=(j-1)(2p+1)+m$, then the
set $A_{i,t}$ will be the main part of the future set $M(1,j,m,t)$.
All $A_{i,t}$ lies in one copy of $K_l$ and so  $G(A_{i,t})$ is connected.
By the definition of the affine plane, if $i'=(j'-1)(2p+1)+m'$
and $i'\neq i$, then the sets $A_{i,t}$ and $A_{i',t'}$ intersect.
Our goal now is to add $2p+1$ vertices to each of $A_{(j-1)(2p+1)+m,t}$ to
provide (b) for $M(1,j,m,t)$ and $M(1,j',m',t')$
{\em only for the same} $j$ and $m$.

Let us fix $j$ and $m$.
For every $t=1,\ldots,p^2$, the set
$M(1,j,m,t)$ will be obtained from $A_{(j-1)(2p+1)+m,t}$ by adding
a $(2p+1)$-element subset of $\cup_{r=1}^{2p+1} Q_m(j,r)$, where
$Q_m(j,r)$ is the copy of $Q_m$ that is contained in $K_l(1,j,r)$.
Every
$t=1,\ldots,p^2$ can be written in the form $t=(a_1-1)p+a_2$,
where $1\leq a_1,a_2\leq p$. So, we include into $M(1,j,m,t)$ the entry
$(a_1,a_2)$ of the square $Q_m(j,m)$. We call this vertex $F(a_1,a_2,j,m)$.
Since $F(a_1,a_2,j,m)$ is in the same
copy of $K_l$ as $A_{(j-1)(2p+1)+m,t}$, it is adjacent to every vertex
in this set. Let $C_a(j,r,m)$ and $R_b(j,r,m)$
denote the $a$th column and the $b$th row of the square
$Q_m(j,r)$, respectively. If $t=(a_1-1)p+a_2$, then our set $M(1,j,m,t)$
will consist of $A_{(j-1)(2p+1)+m,t}$, the  vertex $F(a_1,a_2,j,m)$,
 the row $R_{a_1}(j,m+a_2,m)$ and
the column $C_{a_2}(j,m-a_1,m)$, where the values $m+a_2$ and $m-a_1$
are calculated modulo $2p+1$. Since $F(a_1,a_2,j,m)$
is adjacent to the $a_2$-s entry of
the row $R_{a_1}(j,m+a_2,m)$ and to the $a_1$-s entry of
the column $C_{a_2}(j,m-a_1,m)$, condition (a) holds. Since the projection on
$Q_m$ of $R_{a_1}(j,m+a_2,m)\cup C_{a_2}(j,m-a_1,m)$ is a cross, (b) also holds.

This finishes the construction. It implies that
the Hadwiger number of $G=K_h\Box K_l$, with $h\geq l$ is at least
$$\left\lfloor\frac{h}{(p-1)(2p+1)/2}\right\rfloor p^2 (p-1)(2p+1)/2=
(h-O(p^2))p^2.$$

By Corollary~\ref{cor1}, $p^2=(1-o(1))\sqrt{l}$.  Hence the result.
\end {proof}

\noindent
The following Theorem is an immediate consequence of Theorem \ref{result2}
\begin{theorem}
\label{cor2}
Let $G_1, G_2$ be any two graphs with $\mr(G_1) = h$, $\mr(G_2) = l$ and
$\mr(G_1) \ge \mr(G_2)$. Then
$\mr(G_1 \cart G_2) = \mr(G_2 \cart G_1) \ge h\sqrt l \left (1 - o(1) \right )$.
\end{theorem}

\subsection{Tightness of the lower bound}

Let $K_n$ and $K_m$ be the complete graphs on $n$ and $m$ vertices respectively ($n \ge m$)
and let $h$ be the maximum number such that $K_h \minor K_n \cart K_m$. Let the sets
$V_0,\cdots,V_{h-1} \subseteq V(K_n \cart K_m)$ be the pre-images of vertices of $K_h$ in
 $K_n \cart K_m$. Thus, the vertex sets $V_0,\cdots,V_{h-1}$ are pairwise  disjoint and pairwise
adjacent. Moreover $V_i$   $(0 \le i \le h-1)$ induces a connected subgraph in $K_n \cart K_m$.

Without loss of generality let  $n_0= |V_0| = \min_{0 \le i \le h-1} { |V_i|} $.
Thus $n_0 \le \frac {nm}{h}$. For $S \subseteq V(K_m \cart K_n)$
let $N(S) = \bigcup_{u \in S} N(u)  - S$. (Here $N(u)$ denote the neighbors
of $u$ in $K_n \cart K_m$.)  Since $K_h$ is a complete graph minor of $K_m \cart K_n$,
we have:

\begin {eqnarray}
\label {firsteqn}
 |N(V_0)|   \ge h - 1
\end {eqnarray}

Since $V_0$ induces a connected graph in $K_n \cart K_m$,
the vertices of $V_0$ can be ordered as $v_1,\cdots,v_{n_0}$ such that for $2 \le j \le n_0$,
$v_j$ is adjacent to at least one of the vertices in $\brac{v_1,\cdots,v_{j-1}}$.
Let us define a sequence of sets $\emptyset=X_0,X_1, \cdots,X_{n_0}=V_0$
by setting $X_j = X_{j-1} \cup \{v_j\}$, for $1 \le j \le n_0$.
Clearly, $|N(X_1)| = n+m-2$.  We claim that $|N(X_j)| \le |N(X_{j-1})| + n -2$,
for $2 \le j \le n_0$.
To see this, recall that  $v_j$ is adjacent to at least one vertex $v_k \in X_{j-1}$.
Clearly, out of the $n+ m -2$ neighbors of $v_j$, at least $m-2$ are neighbors of
$v_k$ also, and thus are already in $N(X_{j-1})$. Now, accounting for $v_j$ and $v_k$ also,
 we have $|N(X_j)| \le |N(X_{j-1})| + n -2$, as required. Thus we get
$|N(V_0)|= |N(X_{n_0})|  \le n+ m-2 + (n_0 - 1) (n-2) \le n+m-2 + (\frac {nm}{2} -1) (n-2)$.
Combining this with Inequality \ref {firsteqn} we get:

$$ n + m -2 + ( \frac {nm}{h} -1)(n-2) \ge  h-1 $$

It is easy to verify that if $h > n \sqrt m +   m$, the above inequality will not be satisfied.
So we infer that $h \le n \sqrt m + m$. (Recalling $n \ge m$,  the upper bound tends to
 $n \sqrt m$ asymptotically.)

\subsection{Nonexistence of an upper bound that depends only on $\mr(G_1)$ and $\mr(G_2)$}

We have seen that if we take $G_1$ and $G_2$ as the complete graphs on $k_1$ and $k_2$ vertices
respectively ($k_1 \ge k_2$), then $\mr(G_1 \cart G_2) \le  k_1\sqrt{k_2} + k_2$ for some constant
$c$. It is very natural to ask the following question. Let $G_1$ and $G_2$ be two arbitrary
graphs with $\mr(G_1) = k_1$ and $\mr(G_2) = k_2$. Then does there exists a function
$f:N \times N \rightarrow N$, such that $\mr(G_1 \cart G_2) \le f(k_1,k_2)$? In this
section we demonstrate that in general such a function cannot exist.

\noindent
\begin{definition}(Grid)
 An $n \times n$ grid is a graph with the vertex
set $V = \brac{1,\cdots,n} \times \brac{1,\cdots,n}$. Nodes $\node{i,j}$ and
$\node{i',j'}$ are adjacent if and only if $|i-i'| + |j-j'| = 1$. Note that, an
$n \times n$ grid (which can be viewed as the adjacency graph on an $n \times n$
chessboard) has $n$ rows and $n$ columns, where $i$th row is the induced path
on the vertex set $\brac{\node{i,1},\cdots,\node{i,n}}$ and $j$th column is the induced
path on the vertex set $\brac{\node{1,j},\cdots,\node{n,j}}$.  \\
\end{definition}

\noindent
\begin{definition}(Double-grid)
 An $n \times n$ double-grid is obtained by taking
two $n \times n$ grids and connecting the identical vertices (vertices with
identical labels) from the two grids by an edge.\\
\end{definition}

\noindent
Let $R_n$ be an $n \times n$ grid.
It is easy to see that $R_n$ is a planar graph and hence
$\mr(R_n) \le 4$. By the definition of Cartesian product, $R_n \cart K_2$ is an
$n \times n$ Double-grid. It was proved in \cite{sunns1} that the Hadwiger
number of an $n \times n$ double-grid is at least $n$.
(We give here a sketch of their proof. Let $G_1$ and $G_2$ be the two grids of
the double grid $R_n \cart K_2$. Observe that there is an edge
between any ``row" of $G_1$ and any ``column" of $G_2$. Contracting all the rows of $G_1$ and
all the columns of $G_2$ we get a complete bipartite graph $K_{n,n}$, from which
we easily obtain a $K_n$ minor)

Thus, $\mr(R_n \cart K_2)\ge n$, while $\mr(R_n) \le 4$ and $\mr(K_2) = 2$. This
example shows that in general there is no upper bound on $\mr(G_1 \cart G_2)$ which
depends only on $\mr(G_1)$ and $\mr(G_2)$.

\subsection{Consequences of  Theorem \ref{cor2}. Hadwiger's conjecture for graph products}

\subsubsection{In terms of chromatic number}

Theorem \ref{cor2} naturally leads us to the following question:
Let $G_1$ and $G_2$ be any two graphs with $\chi(G_1) = k_1$ and
$\chi(G_2) = k_2$, where $k_1 \ge k_2$. Let $f(k_1)$ be such that if
$k_2 \ge f(k_1)$, Hadwiger's conjecture is true for $G_1 \cart G_2$. In
fact Hadwiger's conjecture states that $f(k_1) = 1$. Since Hadwiger's
conjecture in the most general case, seems to be hard to prove, it is
interesting to explore how small we can make $f(k_1)$, so that the conjecture
can still be verified, for $G_1 \cart G_2$. To obtain a bound on $f(k_1)$, we
need the following result, proved by Kostochka \cite{Kostochka82} and
 Thomason \cite{AGThoma}, independently.

\begin{lemma}
\label{Kos}
For any graph $G$, $\mr(G) \ge \frac{c_2\chi(G)}{\sqrt{\log{\chi(G)}}}$, where
$c_2$ is a constant.
\end{lemma}

\noindent
As a consequence of Theorem \ref{cor2}, we have the following result.

\begin{theorem}
\label{remark1}
Let $G_1$ and $G_2$ be any two graphs.There exists a constant $c'$ such that if
 $\chi(G_1) \ge \chi(G_2) \ge c' {{\log}^{1.5}}(\chi(G_1))$,
then Hadwiger's conjecture is true for $G_1 \cart G_2$.
\end{theorem}

\begin{proof}
Let $k_1 = \chi(G_1)$ and $k_2  = \chi(G_2)$. Applying Lemma \ref{Kos}
and Theorem \ref{cor2} and noting that $\sqrt{\sqrt{\log(k_2)}} \le
(\sqrt{\log(k_1)})^{0.5}$, we have

$$\mr(G_1 \cart G_2) \ge {c_1}{c_2}^{1.5} \frac{k_1\sqrt{k_2}}
{(\sqrt{\log{k_1}})^{1.5}}$$
Now taking $c' = \frac{1}{(c_1 {c_2}^{1.5})^2}$, ($c_1$ and $c_2$ are the
constants that correspond to  Theorem \ref{cor2} \footnote { From Theorem  \ref {cor2}
we have $\mr (G_1 \cart G_2) \ge c_1 h \sqrt l$, where $c_1$ is a constant.}  and Lemma \ref{Kos}
respectively) and recalling that
$k_2 \ge c'{{\log}^{1.5}}(k_1)$, we get $\mr(G_1 \cart G_2)$ $\ge k_1
= \chi(G_1 \cart G_2)$. The latter equality follows from Lemma \ref{max}.
\end{proof}

\subsubsection{In terms of product dimension}

Recall that the product dimension of a connected graph $G$ is the number of prime
factors in its (unique) prime factorization. It was shown in \cite{sunns1} that if
the product dimension of $G$ is at least $2\log{\chi(G)} +3$, then Hadwiger's
conjecture is satisfied for $G$. Using theorem \ref{cor2}, we can bring this bound
to $2\log{\log{\chi(G)}} +c'$, where $c'$ is a constant. The following Lemma
proved in \cite{sunns1} gives a lower bound for the Hadwiger number of the
$d$-dimensional Hypercube, $H_d$.

\begin{lemma}
\label{lowerbound}
$\mr(H_k) \ge 2^{\floor{(k-1)/2}} \ge 2^{(k-2)/2}$
\end{lemma}

\begin{theorem}
\label{remark2}
Let $G$ be a connected graph and let the (unique) prime factorization of $G$ be
$G = G_1 \cart G_2 \cart ... \cart G_k$. Then there exists a constant $c'$, such
that Hadwiger's conjecture is true for $G$, if $k \ge 2\log{\log{\chi(G)}}
+c'$
\end{theorem}

\begin{proof}
Let $c' = 4\log{\frac{1}{c_1c_2}}+3$, where $c_1$ and $c_2$ are the constants
\footnote{$c_1 ,c_2 \le 1$. So $\frac{1}{c_1c_2} \ge 1$.}
from Theorem \ref{cor2} and Lemma \ref{Kos} respectively.

\noindent
We may assume that $\chi(G_1) \ge \chi(G_{i})$,
for all $ i > 1$.
By Lemma \ref{max}, $\chi(G) =
\max\brac{\chi(G_1), \chi(G_2), ..., \chi(G_k)} = \chi(G_1)$.

\noindent
Let $X = G_2 \cart G_3 \cart \cdots \cart G_{k}$.
Since $G$ is connected, each $G_i$ is also connected. Moreover, $G_i$ has
at least two vertices (and hence at least one edge) since $G_i$ is
prime. It follows that the $(k-1)$-dimensional hypercube is a minor of
$X$.  Thus by Lemma \ref{lowerbound},
$\mr(X) \ge \mr(H_{k-1}) \ge 2^{(k-3)/2} \ge 2^{\log{\log{\chi(G)}} + 2\log{\frac{1}{c_1c_2}}}$.\\

\noindent
Applying Theorem \ref{cor2} to $G_1 \cart X$, we get
$$\mr(G) = \mr(G_1 \cart X) \ge c_1\mr(G_1)\sqrt{\mr(X)}$$
\noindent
Recalling that (by Lemma \ref{Kos}),
$\mr(G_1) \ge \frac{c_2\chi(G_1)}{\sqrt{\log{\chi(G_1)}}}$, we get $\mr(G) \ge \chi(G)$.
\end{proof}

\section{Hadwiger's Conjecture for $G_1 \cart G_2$ when $\chi(G_1) = \chi(G_2)$}

Theorem \ref{remark1} implies the following. Let $G_1$ and $G_2$ be two graphs
such that $\chi(G_1) = \chi(G_2)$. Then $G_1 \cart G_2$ satisfies Hadwiger's
conjecture if $\chi(G_1) = \chi(G_2) = t$ is sufficiently large.($t$ has to be
sufficiently large, because of the constant $c'$ involved in Theorem \ref{remark1}).
In this section we give a different proof for this special case.
We show that irrespective of the value of $t$ ($ = \chi(G_1)$), $G_1 \cart G_2$
satisfies Hadwiger's conjecture if $\chi(G_1) = \chi(G_2)$.

A graph $G$ is said to be $k$-critical if and only if
$\chi(H) < \chi(G)$ for every proper subgraph $H$ of $G$.
Every $k$-chromatic graph has a $k$-critical subgraph in it, obtained by
greedily removing as many vertices and edges as possible from $G$, such that the
chromatic number of the resulting graph remains the same.

\noindent
We need the following two Lemmas, the proofs of which can be
found in \cite{West}.

\begin{lemma}
\label{degree}
If $G$ is a k-critical graph, then the minimum degree of $G$, $\delta(G) \ge k-1$
\end{lemma}

\begin{lemma}
\label{Dirac}
Let $G$ be a graph with minimum degree $\delta$.
Then $G$ contains a simple path on at least $\delta+1$ vertices.
\end{lemma}

\noindent
We use $W_n$ to denote the graph whose vertex set is $\brac{0,1,...,n-1}$ with
an edge defined between two vertices $i$ and $j$ (assuming $i < j$) if and only if
either $i = 0$ or $j=i+1$.
$W_n$ is essentially a simple path on $n$
vertices, with the extra property that vertex $0$ is adjacent to all the other
vertices. An illustration of $W_n$ is given in Figure 1.

\begin{center}
\begin{figure}[h]
\scalebox{1.00}{\includegraphics{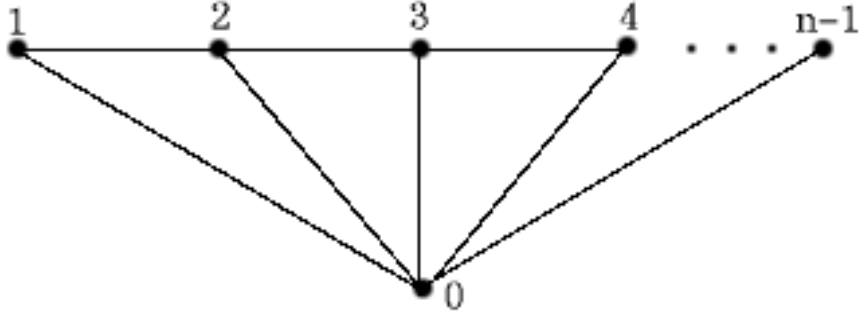}}
\caption{Illustration of $W_n$}
\end{figure}
\end{center}

\begin{lemma}
\label{wk}
Every $k$-chromatic graph $G$ has $W_k$ as a minor.
\end{lemma}

\begin{proof} Let $H$ be a $k$-critical subgraph of $G$.  By Lemma \ref{degree}
$\delta(H) \ge k-1$. Let $P=(v_0,v_1,\cdots,v_{l-1})$
be the longest simple path in $H$. By Lemma \ref{Dirac}, $l \ge k=\chi(G)$.
Let $N(v_0)$ denote the set of neighbors of $v_0$ in $H$. i.e.,
$N(v_0) = \brac{u \in V(H)-{v_0}: (u,v_0) \in E(H)}$.
Since $P$ is the longest simple path, $N(v_0) \subseteq V(P)-\brac{v_0}
=\brac{v_1,v_2,\cdots,v_{l-1}}$. Otherwise if $w \in N(v_0)$ and $w \not\in
V(P)-\brac{v_0}$, then $(w,v_0,v_1,\cdots,v_{l-1})$ will be a longer simple
path in $G$, contradicting the assumption
that $P$ is the longest. Let $\brac{v_{i_1}, v_{i_2},...,v_{i_{k-1}}}
\subseteq V(P)$ be any $k-1$ neighbors of $v_0$ in $H$, where $i_1 \leq
i_2 \leq ... \leq i_{k-1}$. Consider the $k-1$ sub-paths of $P$,
from $v_0$ to $v_{i_1}$, from $v_{i_1}$ to $v_{i_2}$, $\cdots$, from
$v_{i_{k-2}}$ to $v_{i_{k-1}}$.
Contracting each sub-path to a single edge, we get $W_k$ as a minor of $G$.
\end{proof}

\noindent
$W_n \cart W_n$ is the graph with vertex set $V = \brac{0,1,...,n-1}
\times \brac{0,1,...,n-1}$. By the definition of graph Cartesian product,
vertices $\node{i,j}$ and $\node{i',j'}$ are adjacent in $W_n \cart W_n$ if
and only if either $i = i'$ and $(j,j') \in E(W_n)$ or $j = j'$ and
$(i,i') \in E(W_n)$. Thus $\node{i,j}$ and $\node{i',j'}$ in $W_n \cart W_n$
are adjacent if and only if at least one of the following conditions hold.
\\

\noindent
(1). $ i = i'$ and $ j = j' \pm 1$  (2). $ i = i'$ and $ j = 0$
(3). $ i = i'$ and $ j'= 0$\\
\noindent
(4). $ j = j'$ and $ i = i' \pm 1$  (5). $ j = j'$ and $ i = 0$
(6). $ j = j'$ and $ i'= 0$\\

\begin{lemma}
\label{knminor}
$K_n \preceq W_n \cart W_n$.
\end{lemma}

\begin{proof}
For $ 0 \le i \le n-1$, let $B_i \subseteq V(W_n \cart W_n)$ be
defined as $B_i = \brac{\node{i,0},
\node{i,1},...,\node{i,i-1}, \node{i,i}, \node{i-1,i},...,\node{1,i},
\node{0,i}}$. The following properties hold for $B_i$.

\begin{enumerate}
\item For $i \not= j$, $B_i \cap B_j = \emptyset$. This follows from the
definition of $B_i$.
\item Each $B_i$ induces a connected graph. This follows from the fact that
$(\node{i,j}, \node{i,j+1}) \in E(W_n \cart W_n)$ and
$(\node{j,i}, \node{j-1,i}) \in E(W_n \cart W_n)$, by the definition of $W_n \cart W_n$.
\item For $i < j$, $B_i$ and $B_j$ are adjacent. This is because,
$\node{i,0} \in B_i$ , $\node{i,j} \in B_j$ and $(\node{i,0}, \node{i,j})
\in E(W_n \cart W_n)$.
\end{enumerate}

\noindent
In other words, the sets $B_i$ are connected, disjoint and are pair-wise
adjacent. Thus contracting each $B_i$ to a single vertex we get a $K_n$
minor.
\end{proof}

\begin{theorem}
\label{result1}
If $ \chi(G) = \chi(H) $, then Hadwiger's conjecture is true
for $G \cart H $.
\end{theorem}

\begin{proof}
Let $\chi(G) = \chi(H) = n$. By Lemma \ref{wk}, we have
$W_n \minor G$ and $W_n \minor H$. Now Lemma \ref{minor} implies
$W_n \cart W_n \minor G \cart H$. Since by Lemma \ref{knminor},
$K_n \minor W_n \cart W_n$, we have $K_n \minor G \cart H$. This
together with Lemma \ref{max}, gives $\mr(G \cart H)$ $\ge n = \chi(G \cart H)$,
proving the Theorem.
\end{proof}

\noindent
It was shown in \cite{sunns1} that if a graph $G$ is isomorphic to $F^d$, for
some graph $F$ and $d \ge 3$ then Hadwiger's conjecture is true for $G$.
The following improvement is an immediate consequence of Theorem \ref{result1}
 and Lemma \ref{max}.

\begin{theorem}
\label{cor11}
Let a graph $G$ be isomorphic to $F^d$ for some graph $F$ and for $d \ge 2$.
Then Hadwiger's conjecture is true for $G$.
\end{theorem}


\begin{thebibliography}{10}

\bibitem{auren}
F.~Aurenhammer, J.~Hagauer, and W.~Imrich.
\newblock Cartesian graph factorization at logarithmic cost per edge.
\newblock {\em Computational Complexity}, 2:331--349, 1992.

\bibitem{sunns1}
L.~Sunil Chandran and Naveen Sivadasan.
\newblock On the {H}ardwiger's conjecture for graph products.
\newblock {\em Discrete Mathematics}, 307(2):266--273, 2007.

\bibitem{Diest}
R. Diestel.
\newblock {\em Graph Theory}, volume 173.
\newblock Springer Verlag, New York, 2 edition, 2000.

\bibitem{DieRem}
R. Diestel and C.~Rempel.
\newblock Dense minors in graphs of large girth.
\newblock {\em Combinatorica}, 25:111--116, 2005.

\bibitem{Dirac1}
G.~A. Dirac.
\newblock In abstrakten {G}raphen vorhandene vollst\"andige $4$--{G}raphen und
  ihre {U}nterteilungen.
\newblock {\em Math. Nachr.}, 22:61--85, 1960.


\bibitem{Hadwiger43}
H.~Hadwiger.
\newblock {\"Uber eine Klassifikation der Streckenkomplexe, Vierteljscr}.
\newblock {\em Naturforsch. Gessellsch. Z\"urich}, 88:132--142, 1943.


\bibitem{Sandy}
W. Imrich and S. Klav\^zar.
\newblock {\em Product Graphs: Structure and Recognition}.
\newblock John Wiley and Sons,Inc, 2000.

\bibitem{IP}
H.~Iwaniec and J.~Pintz.
\newblock Primes in short intervals.
\newblock {\em Monatsh. Math.}, 98:115--143, 1984.

\bibitem{Kostochka82}
A.~V. Kostochka.
\newblock The minimum hadwiger number of graphs with a given mean degree of
  vertices.
\newblock {\em Metody Diskret. Analiz.}, 38:37--58, 1982.
\newblock {(In Russian)}.

\bibitem{Kotlov}
A. Kotlov.
\newblock Minors and strong products.
\newblock {\em European Journal of Combinatorics}, 22:511--512, 2001.

\bibitem{KuhOs}
D. K\"uhn and D. Osthus.
\newblock Minors in graphs of large girth.
\newblock {\em Random Structures and Algorithms}, 22:213--225, 2003.

\bibitem{Mad2}
W.~Mader.
\newblock Homomorphies\"atze f\"ur {G}raphen.
\newblock {\em Math. Annalen}, 178:154--168, 1968.

\bibitem{Miller1}
Z.~Miller.
\newblock Contractions of graphs: A theorem of {O}re and an extremal problem.
\newblock {\em Discrete Mathematics}, 21:261--273, 1978.


\bibitem{RobertsonSeymourThomas93}
N.~Robertson, P.~D. Seymour, and R.~Thomas.
\newblock Hadwiger's conjecture for {K6-free} graphs.
\newblock {\em Combinatorica}, 13:279--361, 1993.

\bibitem{R}
H.~J. Ryser.
\newblock {\em Combinatorial mathematics,The Carus Mathematical Monographs},
  volume~14.
\newblock Wiley, New York, 1963.

\bibitem{Sabi57}
G.~Sabidussi.
\newblock Graphs with given group and given graph theoretic properties.
\newblock {\em Canad. J. Math.}, 9:515--525, 1957.

\bibitem{AGThoma}
A.~G. Thomason.
\newblock An extremal function for contractions of graphs.
\newblock {\em Math. Proc. Camb. Phil. Soc.}, 95:261--265, 1984.



\bibitem{Wagner2}
K.~Wagner.
\newblock {\"Uber} eine {Eigenschaft} der {Ebenen} {Komplexe}.
\newblock {\em Math. Ann.}, 114:570--590, 1937.

\bibitem{Wagner1}
K.~Wagner.
\newblock Beweis einer Abschw\"achung der Hadwiger--Vermutung.
\newblock {\em Math. Annalen}, 153:139--141, 1964.

\bibitem{West}
D.~B. West.
\newblock {\em Introduction to Graph Theory}.
\newblock Prentice Hall India, NewDelhi, 2 edition, 2003.

\end{thebibliography}
\end{document}